\documentclass[12pt,leqno]{article}
\usepackage{ geometry,amsmath, amsfonts}     
\usepackage{latexsym}
\usepackage{mathrsfs,color}           
\newtheorem{theorem}{Theorem}
\newtheorem{proposition}{Proposition}[section]
\newtheorem{lemma}[proposition]{Lemma}
\newtheorem{corollary}[proposition]{Corollary}

\newtheorem{remark}[proposition]{Remark}

\numberwithin{equation}{section}
\newcommand{\bpf} {\noindent{\sc Proof} : }
\newcommand{\epf}  {\hfill$\diamondsuit$\vspace{.5cm}}
\newcommand{\E} {\mathbb{E}}
\renewcommand{\P} {\mathbb{P}}

\newcommand{\Q}{\mathbb{Q}}
\newcommand{\R}{\mathbb{R}}
\newcommand{\Z}{\mathbb{Z}}
\newcommand{\F} {\mathcal{F}}
\newcommand{\M} {\mathcal{M}}
\newcommand {\ep} {\varepsilon}
\newcommand{\eps} {\varepsilon}
\newcommand {\x} {\mathcal{X}}

\title{Muller's ratchet clicks in finite time}
\author{Julien Audiffren\thanks{Aix-Marseille Universit\'e, LATP 39, rue F. Joliot Curie 13453 Marseille cedex 13, julien.audiffren@wanadoo.fr, partially supported by the ANR project MANEGE}, Etienne Pardoux\thanks{Corresponding author,
Aix-Marseille Universit\'e, LATP 39, rue F. Joliot Curie 13453 Marseille cedex 13, 
Tel~: +33(0)413 55 14 57, Fax~: +33(0)413 55 13 66, pardoux@cmi.univ-mrs.fr, partially supported by the ANR project MANEGE}}
\date{}
\begin{document}
\maketitle
\begin{abstract}
We consider the accumulation of deleterious mutations in an asexual population, a phenomenon known as Muller's ratchet, using the continuous time model proposed in \cite{epw}. We show that for any parameter $\lambda>0$ (the rate at which mutations occur), for any $\alpha>0$ (the toxicity of the mutations) and for any size $N>0$ of the population, the ratchet clicks a.s. in finite time. That is to say the minimum number of deleterious mutations  in the population goes to infinity a.s.
\end{abstract}

\section{Introduction}

In natural evolution, deleterious mutations occur much more frequently than beneficial ones. Since the last category is always favored by selection, one may wonder about the advantage of sexual reproduction over the asexual type. The answer has been proposed : in an asexually reproducing population, each individual always inherits all the deleterious mutations of his ancestor (except if another mutation occurs at the same locus on the genome; but this event is rare and we will not consider it),  whereas in sexual reproduction, recombinations occur, which allow an individual to take part of a chromosome from each of his parents, therefore giving him a chance to get rid of deleterious mutations. Muller's ratchet can be used as an attempt to translate this phenomenon in a mathematical model, thus explaining the advantage of sexual reproduction \cite{may}. If one considers the best class (the group of fittest individuals)  in a given asexual population, Muller's ratchet is said to click when the best class becomes empty. Since beneficial mutations do not occur in this model, it means that all the individuals of the best class  have mutated. 

The first model for Muller's ratchet  due to Haigh \cite{hai} can be described as follows. Consider an asexual population of fixed sized $N$ which evolves in discrete time, with a multiplicative selection model. Only deleterious mutations occur. Denoting by $0 \le\alpha \le 1$ the deleterious strength of the mutations, and by $\lambda >0$ the rate at which they occur, every generation is constituted as follows : each individual chooses a parent from the previous generation, in such a way that the probability of choosing a specific father with $k$ deleterious mutations is (we denote by $N_k$ the number of such individuals in the previous generation) :
$$ \frac{ (1-\alpha)^{k}}{\sum_{k=0}^\infty N_k(1-\alpha)^k}.$$ Next each newborn gains $\xi$ deleterious mutations, where $\xi$
is a Poisson random variable with parameter $\lambda$. It is immediate to see that  this model clicks a.s. in finite time. Indeed at each generation, with probability $(1-\exp(-{\lambda}))^N$ all the individuals mutate, which induces the click. 

There are three parameters in our model~:

$N $ is the size of the population,

$\lambda$ is the mutation rate,     

$\alpha$ is the fitness decrease due to each mutation.

The Fleming--Viot model for Muller's ratchet proposed by A. Etheridge, P. Pfaffelhuber and A. Wakolbinger in \cite{epw} consists of the following infinite set of SDEs for the 
$X_k(t)$'s, $k\ge0$, where $X_k(t)$ denotes the proportion of individuals in the population 
who carry exactly $k$ deleterious mutations at time $t$ (with $X_{-1}\equiv0$) :
\begin{equation}\label{un}
\left\{ 
\begin{aligned}
\!\! dX_k(t)  &=\!  \left[  \alpha(M_1 (\!t\!)\! -\! k)X_k (\!t\!) + \lambda (X_{k-1} (\!t\!)\! -\! X_k (\!t\!)) \right]\! dt\! +\!\!\!\!\! \sum_{\ell \geq 0, \ell \neq k }\!\!\!\! \sqrt{ \frac{X_k(t) X_\ell(t) } {N} } dB_{k,\ell} (\!t\!), \\
\!\! X_k(0)&=x_k,\ k\ge 0;
\end{aligned}
\right.
\end{equation}
where $\left\{ B_{k,\ell}, k > \ell \ge 0 \right\}$ are independent Brownian motions,  $B_{k,\ell}=-B_{\ell,k}$; and $M_1(t)$ $=\sum_{k\ge0}kX_k(t)$.

The first term in the drift models the selective effect of the deleterious mutations. Those individuals
who carry less (resp. more) mutations than the average number of mutations in the population have a selective advantage  (resp. disadvantage). The second term in the drift reflects the effect of the accumulation of mutations : at rate $\lambda$, individuals carrying $k-1$ mutations gain a $k$--th mutation, they jump into the $k$--class, and at the same rate individuals carrying $k$ mutation gain a $k+1$--th mutation, they jump out of the $k$--class. The diffusion term reflects the resampling effect of the  birth events, where the factor $N^{-1/2}$ can be understood as being equivalent to the rescaling of time $t\to t/N$, if $N$ is the ``effective population size'', which is natural in Kingman's 
coalescent \cite{kin}. For the equivalence between the present model and a more intuitive look--down model \`a la Donnelly--Kurtz, we refer the reader to \cite{aud}. 

We will show in section \ref{sec2} that the infinite dimensional system of SDEs \eqref{un} is well posed provided we choose the initial condition $x=(x_k,\ k\ge0)\in\x_\delta$ for some $\delta>0$, where
\begin{equation}\label{DefSet}
\x_\delta:=\left\{x\in[0,1]^\infty,\quad \sum_{k=0}^\infty x_k=1,\quad \sum_{k=0}^\infty k^{2+\delta}x_k<\infty
\right\}.
\end{equation}

\bigskip

We define $ T_0= \inf\{t>0, X_0(t)=0\}$.
The purpose of the present work is to show that this model of Muller's ratchet is bound to click in finite time, that is to say $T_0<\infty $ a.s.  We are going to prove the following theorem :

\begin{theorem}\label{th}
 For any $\delta>0$,
for any choice of initial condition in $\x_\delta$, let $(X_k(t))_{k \in {\Z_+}}$ be the solution of \eqref{un}. Then $\P(T_0<\infty)=1 $.
\end{theorem}
 
 We will in fact prove a stronger result, namely
 \begin{theorem}\label{thdeu} For any $\delta>0$,
for any choice of initial condition in $\x_\delta$, let $(X_k(t))_{k \in {\Z_+}}$ be the solution of \eqref{un}. Then there exists $\overline{\rho}>0$, which depends upon the parameters $N$, $\alpha$ and $\lambda$, such that $\E\left[\exp(\rho T_0)\right]<\infty$, for
all $0<\rho<\overline{\rho}$.
\end{theorem}

Clearly, a model for Muller's ratchet must have the property that the ratchet clicks in finite time. In a sense our result says that the Etheridge--Pfaffelhuber--Wakolbinger model for Muller's ratchet is a reasonable model, in the sense that it exhibits a.s. clicking, as the computer simulations had already shown, see \cite{epw}. Note that once the zero class is empty, the 1--class takes its place, and some time later a second click happens, at which time both the zero class and the 1--class become empty, and so on. Of course, we would like to know more about the time it takes for the ratchet to click. Here we show that it has an exponential moment of some order. We hope to get more precise information in some future work.
 
There are several difficulties in this model. First, it is an infinite system of SDEs which cannot be reduced to a finite dimensional system. Only $X_0$ and $M_1$ enter the coefficients of the equation for $X_0$, but the equation for $M_1(t)$ brings in the second centered moment $M_2(t)=\sum_{k=0}^\infty(k-M_1(t))^2X_k(t)$. The system of SDEs for the centered moments of all orders is infinite as well, the moments of order up to $\ell=2k$ enter the coefficients of the equation for
the $k$--th centered moment, and there is no known solution to it (except in the deterministic case $N=+ \infty$, which is solved in \cite{epw}). In addition, one has $d \left< X_0,M_1 \right>_t = - \frac{M_1(t) X_0(t)}{N}dt$. There is no simple relation between $X_0$ and $M_1$, except that $X_0+M_1 \geq 1$, and $(X_0=1)$ $\Rightarrow$ $( M_1=0)$. But we could have $X_0 \rightarrow 0$ and $M_1 \rightarrow \infty$.  Last but not least, the diffusion coefficient in $dX_k$ is not a Lipschitz function of $X_k$ at $0$ and $1$, and it vanishes at those two points.

In order to prove the theorem, we will use a three--step proof. First, in section 3 we will show that $M_1$ cannot grow too fast with a good probability, 
and we will deduce that for a specific set of initial conditions, the ratchet does click with a strictly positive probability $p_{fin}$, in a given interval of time.

Next, we show in section 4 that the product $X_0 M_1^2$ is bound to come back under $\frac{2 (\lambda +1)}{ \alpha}$ after any time, and we  use all the previous results to deduce that $M_1$ is also bound to return under $\beta= \frac{\lambda }{\alpha}$ after any time, as long as  the ratchet does not click.

Finally in section 5 we  prove that each time $M_1$ gets below $\beta$, the ratchet clicks with a positive probability in a prescribed interval of time. We then conclude with the help of the strong Markov property.

In section 6 we show how the proof of Theorem \ref{th} can be modified into a proof of Theorem \ref{thdeu}. The reader may wonder why we do not prove Theorem \ref{thdeu} from the very beginning, and first prove a weaker result. The reason is that 
the difference between the two proofs is essentially that while proving Theorem \ref{th}, we prove that as long as the ratchet has not clicked, $M_1$ is bound to return below the value $\beta$, i.e. the drift of $X_0$ is bound to become non--positive, which is an interesting result in itself, while the proof of Theorem \ref{thdeu} is based on the same strategy, but with $\beta$ replaced by a much less explicit quantity.

We shall essentially work with the two dimensional process $ \{ X_0(t), M_1(t) \}$, and we shall use the equation for $X_1$ only in one place, namely in Lemma \ref{X01} in order to show that $X_0$ does not get stuck near the value 1. This means that we shall make use only of the three following equations.
\begin{equation}\label{3eq}
\left\{
\begin{aligned}
\!\! dX_0(t)&=  \left(  \alpha M_1(t) - \lambda  \right) X_0(t) dt + \sqrt{ \frac{ X_0(t) (1- X_0(t)) } {N} } dB_{0}(t),\\
\!\! dX_1(t)&=\!\left(  \alpha( M_1(t)\! -\!1 )X_1(t)\! +\! \lambda( X_0(t)\!-\!  X_1(t)) \right)\! dt\! \!+ \!\!\!
 \sqrt{\! \frac{ X_1(t) (1- X_1(t)) } {N} } dB_{1}\!(\!t\!),\\
\!\! dM_1(t)&=(\lambda-\alpha M_2(t))dt+\sqrt{\frac{M_2(t)}{N}}dB(t).
\end{aligned}
\right.
\end{equation}
The three Brownian motions $B_0$, $B_1$ and $B$ are standard Brownian motions. They
are not independent, and the three dimensional process $(B_0(t), B_1(t), B(t))$ is not a Gaussian process. But this will play no role in our analysis.
This system is not closed, since $M_2$ enters the coefficients of the last equation. However, the crucial remark is that 
it will not be necessary to estimate $M_2$, in order to estimate $M_1$. This is due to the fact that the $M_1$--equation takes the form
$dM_1(t)=\lambda dt+dZ_t$, where $Z_t=W(A_t)-\alpha N A_t$, if $A_t:=N^{-1}\int_0^t M_2(s) ds$ and $\{W(t),\ t\ge0\}$ is a standard Brownian motion. The larger $M_2$ is, the more likely $Z_t$ is negative, which produces a smaller $M_1$. This means that we should be able to estimate $M_1$, without having to estimate $M_2$, which is done below in  Lemma \ref{M_1 pas trop rapide} and \ref{return}. In particular, we show in Lemma \ref{Mini} below that, as long as the ratchet has not clicked, $M_1$ is bound to return below the level $\beta=\lambda/\alpha$ after any time. We believe that this is an interesting qualitative property of the model. Note that Theorem \ref{thdeu} is proved by essentially the same argument as Theorem \ref{th}, but with that level 
$\beta$ replaced by $2\beta\vee(\eps/\delta)$, where the constants $\eps$ and $\delta$, which are defined in the proof of Theorem \ref{th}, have no explicit relation to the constants of the model.

\section{Preliminary results}\label{sec2}

The aim of this section is to establish a weak existence and uniqueness result for the infinite system of SDEs \eqref{un}, under the  condition that the initial condition $\{X_k(0),\, k\ge0\}$ belongs to the set $\x_\delta$ for some $\delta>0$ (see \ref{DefSet} for the definition of this set).

We equip this set with the topology under which a probability $x^n=(x^n_k,\ k\ge0)$ on
$\mathbb{Z}_+$ converges to
$x=(x_k,\ k\ge0)$ if both it converges weakly, and $\sup_n\sum_{k\ge0}k^{2+\delta}x^n_k<\infty$.
More precisely, we will prove in this section
\begin{theorem}\label{exist-uniq}
If the initial condition $x$ belongs to $\x_\delta$, for some $\delta>0$, then \eqref{un} has a unique weak solution $X(t)=\{X_k(t),\ k\ge0\}$ which is a. s. continuous with values in $\x_\delta$.
\end{theorem}
\begin{remark}
Previous results on this system of SDEs assume that the probability $x$ on
$\Z_+$ possesses exponential moments of arbitrary order, see \cite{cuth}, or of some
order, see \cite{psw}. This assumption is naturally requested if one wants to be able to
write equations for arbitrary moments of the random measure $X(t)$ on $\Z_+$. However,
we will need only to make sure that $M_1(t)$ and $M_2(t)$  have finite expectation, and for that purpose our weaker condition will be sufficient.
\end{remark}

We start with the case $\alpha=0$.
\subsection{The case $\alpha=0$}
\begin{proposition}\label{alpha=0}
Suppose that $\alpha=0$. Then, for any
 initial condition $x\in\x_\delta$, \eqref{un} has a unique weak solution $X(t)=\{X_k(t),\ k\ge0\}$ which is a. s.  continuous with values in $\x_\delta$, and is such that for each $\lambda, \delta>0$, there exists a
 locally bounded function $C_{\lambda,\delta}(t)$ such that
 \begin{equation}\label{deltamoment}
 \E\sum_{k=0}^\infty k^{2+\delta}X_k(t)\le C_{\lambda,\delta}(t).
 \end{equation}
\end{proposition}
\bpf
Let us rewrite our system of SDEs in the particular case $\alpha=0$ (again it is written with the convention that $X_{-1}(t)\equiv0$) in the form
\begin{equation}\label{un0}
\left\{ 
\begin{aligned}
dX_k(t)  &=   \lambda (X_{k-1}(t) - X_k(t))  dt +  d\M_k(t),\ k\ge0; \\
\langle\M_k,\M_\ell\rangle_t&=N^{-1}\int_0^t X_k(s)(\delta_{k,\ell} - X_\ell(s))ds,\ k,\ell\ge0;\\
 X_k(0)&=x_k,\ k\ge 0;
\end{aligned}
\right.
\end{equation}
where the $\M_k(t)$'s are continuous martingales, and $\langle\M_k,\M_\ell\rangle$ stands for the joint quadratic variation of the two martingales $\M_k$ and $\M_\ell$.
We can apply the result of Theorem 2.1 in \cite{shi}, which ensures that \eqref{un0}
has a unique weak solution. The facts that
 $X_k(t)\ge0$, for all $k\ge0$, $t\ge0$, a.s. and $\sum_{k\ge0}X_k(t)=1$ for all $t\ge0$ a.s.
 follow from the results in \cite{shi}. 
 
 We now have
\begin{align*}
\E\left(\sum_{k=0}^K kX_k(t)\right)&=\sum_{k=0}^K kx_k+
\lambda\E\int_0^t\sum_{k=0}^K \left(kX_{k-1}(s)-kX_k(s)\right)ds,\\
\E\left(\sum_{k=0}^\infty kX_k(t)\right)&\le \sum_{k=0}^\infty kx_k+\lambda t,
\end{align*}
since $\sum_{j=0}^{K-1}X_j(s)\le1$. 
Furthermore, using this last inequality in the last step below,
\begin{align*}
\E\left(\sum_{k=0}^K k^2X_k(t)\right)&=\sum_{k=0}^K k^2x_k+
\lambda\E\int_0^t\sum_{k=0}^K \left(k^2X_{k-1}(s)-k^2X_k(s)\right)ds\\
&\le \sum_{k=0}^K k^2x_k+\lambda\E\int_0^t\sum_{j=0}^{K-1}(2j+1)X_j(s)ds,\\
\E\left(\sum_{k=0}^\infty k^2X_k(t)\right)&\le \sum_{k=0}^\infty k^2x_k+\lambda t+\lambda^2 t^2+2\lambda t\sum_{k=0}^\infty kx_k.
\end{align*}
Let us now suppose that $0<\delta\le1$, and we exploit the fact that 
${2+\delta}k^{1+\delta}\le 3k^2$. We then deduce that
\begin{align*}
\E\left(\sum_{k=0}^K k^{2+\delta}X_k(t)\right)&=\sum_{k=0}^K k^{2+\delta}x_k+
\lambda\E\int_0^t\sum_{k=1}^K \left(k^{2+\delta}X_{k-1}(s)-k^{2+\delta}X_k(s)\right)ds\\
&\le \sum_{k=0}^K k^{2+\delta}x_k+3\lambda\E\int_0^t\sum_{k=1}^{k+1}j^2X_k(s)ds,\\
\E\left(\sum_{k=0}^\infty k^{2+\delta}X_k(t)\right)&\le  C_2(\lambda,t),
\end{align*}
from the last estimate. If $\delta>1$, we need to estimate the third moment in terms of the second,
then the fourth in terms of the third, ..., and finally the $2+\delta$--th in terms of the 
$2+\lfloor\delta\rfloor$--th.

So far we have proved that $X(t)\in\x_\delta$ a. s. for all $t\ge0$. We now prove that in fact
a. s., $X(t)\in\x_\delta$ for all $t\ge0$. Our next argument will be very similar to an argument in \cite{psw}. For any $m\ge1$, $t\ge0$, let
$$N_{m,\delta}(t):=\sum_{k=0}^\infty \inf(k,m)^{2+\delta}X_k(t).$$
It is easy to check that $\{N_{m,\delta}(t),\ t\ge0\}$ is a positive submartingale, to which we
can apply Doob's inequality, which, together with the monotone convergence theorem, yields
that for any $K, T >0$, 
\begin{align*}
\P\left(\sup_{0\le t\le T}\sum_{k=0}^\infty k^{2+\delta}X_k(t)>K\right)
&=\lim_{m\to\infty}\P\left(\sup_{0\le t\le T}N_{m,\delta}(t)>K\right)\\
&\le\lim_{m\to\infty}K^{-1}\E\left[N_{m,\delta}(T)\right]\\
&=K^{-1}\E\left[\sum_{k=0}^\infty k^{2+\delta}X_k(T)\right]\\
&\le K^{-1}C_{\lambda,\delta}(T),
\end{align*}
where we have used \eqref{deltamoment} for the last inequality. It now follows that for all $T>0$,
$$ \P\left(\sup_{0\le t\le T}\sum_{k=0}^\infty k^{2+\delta}X_k(t)<\infty\right)=1.$$
The a. s. continuity with values in $\x_\delta$ is now easy to check.
\epf 
 
 We next want to establish the equation for the first moment $M_1(t):=\sum_{k\ge1}kX_k(t)$. This equation will involve the process $M_2(t)=\sum_{k\ge1}k^2X_k(t)-[M_1(t)]^2$. We know by now that those quantities are well defined and finite.
 \begin{proposition}\label{prop-M1eq}
 The first moment solves the SDE
 $$dM_1(t)=\lambda dt+d\M(t),$$
 where $\{\M(t),\ t\ge0\}$ is a continuous martingale satisfying
 $$\langle\M,\M\rangle_t=N^{-1}\int_0^tM_2(s)ds,$$ and  for any $k\ge0$,
 \begin{equation}\label{covar}
 \langle \M,\M_k\rangle_t=N^{-1}\int_0^t(k-M_1(s))X_k(s)ds.
 \end{equation}
\end{proposition}
\bpf
For any $K>1$, let $M_{1,K}(t):=\sum_{k=1}^K kX_k(t)$. We have readily
$$M_{1,K}(t)=M_{1,K}(0)+\lambda\int_0^t\sum_{j=0}^{K-1}X_j(s)ds-\lambda\int_0^tKX_K(s)ds
+{\mathcal M}_{1,K}(t),$$
where ${\mathcal M}_{1,K}(t)$ is a continuous martingale, with
\begin{align*}
d\langle {\mathcal M}_{1,K}\rangle_t=N^{-1}\left(\sum_{k=1}^K k^2X_k(t)-\left[M_{1,K}(t)\right]^2\right)dt.
\end{align*}
It follows from \eqref{deltamoment} that 
$$\E\int_0^tKX_K(s)ds\to0$$ as $K\to\infty$. Consequently, all terms in the above equation converge
as $K\to\infty$, yielding that
$$M_1(t)=M_1(0)+\lambda t+{\mathcal M}(t),$$ where ${\mathcal M}(t)$ is a continuous martingale as follows from the next lemma, which is such that
$$d\langle {\mathcal M},\M\rangle_t=N^{-1}M_2(t) dt.$$
Moreover, if $1\le k\le K$,
\begin{align*}
d\langle \M_{1,K},\M_k\rangle_t&=N^{-1}[kX_k(t)(1-X_k(t))-\sum_{\ell\not=k, \ell\le K}\ell X_k(t)X_\ell(t)]dt\\
&=N^{-1}X_k(t)[k-\sum_{\ell\le K}\ell X_\ell(t)]dt.
\end{align*}
The second part of the result follows, by letting $K\to\infty$. \epf

To complete this last proof, we need to establish
\begin{lemma}
The collection of processes $\{\M_{1,K}(t),\ t\ge0\}_{K\ge1}$ is tight in 
$C([0,+\infty))$.
\end{lemma}
\bpf
From the Corollary of Theorem 7.4 page 83 in \cite{bil}, Chebychef's and Doob's inequalities, it suffices to prove that for each $T>0$ there exists a constant $C(\delta,T)$ such that
for all $0\le s<t\le T$,
\begin{equation}\label{tight}
\E\left[|\M_{1,K}(t)-\M_{1,K}(s)|^{2+\delta}\right]\le C(\delta,T)|t-s|^{1+\delta/2}.
\end{equation}
From the well--known Davis--Burkholder--Gundy inequality (see e.g. p. 160 in \cite{ry}), there exists a constant
$c(\delta)$ such that
$$
\E\left[|\M_{1,K}(t)-\M_{1,K}(s)|^{2+\delta}\right]\le c(\delta)\E\left[\left(\langle\M_{1,K}\rangle_t-
\langle\M_{1,K}\rangle_s\right)^{1+\delta/2}\right].
$$
We have, using Jensen's inequality in two distinct instances,
\begin{align*}
\E\left[\left(\int_s^t\sum_{k\ge0}k^2X_k(r)dr \right)^{1+\delta/2}\right]
&\le (t-s)^{\delta/2}\E\int_s^t\left(\sum_{k\ge0}k^2X_k(r)\right)^{1+\delta/2}dr\\
&\le (t-s)^{\delta/2}\E\int_s^t\sum_{k\ge0}k^{2+\delta}X_k(r)dr.
\end{align*}
\eqref{tight} follows by combining the two last estimates with \eqref{deltamoment}. \epf

\subsection{The general case}
We can now prove Theorem \ref{exist-uniq}. We first proceed with the

\noindent{\sc Proof of existence}
We now introduce a Girsanov transformation. It follows from Proposition \ref{prop-M1eq} that there exists a Brownian motion $\{B(t),\ t\ge0\}$ such that
$$dM_1(t)=\lambda dt+\sqrt{\frac{M_2(t)}{N}}dB_t.$$
For any $\alpha>0$, let
$$Z_\alpha(t):=\exp\left(-\alpha\sqrt{N}\int_0^t\sqrt{M_2(s)}dB_s-\frac{\alpha^2N}{2}
\int_0^tM_2(s)ds\right).$$
It is easily seen that
\begin{align*}
Z_\alpha(t)&=\exp\left(N\alpha\left[M_1(0)+\lambda t-M_1(t)-\frac{\alpha}{2}\int_0^tM_2(s)ds\right]\right)\\
&\le\exp\left(N\alpha\left[M_1(0)+\lambda t\right]\right).
\end{align*}
It is now clear that $\{Z_\alpha(t),\ t\ge0\}$ is a martingale, and consequently there exists a unique probability
measure $\P^\alpha$ on $(\Omega,\mathcal{F})$, such that for all $t>0$,
$$\frac{d\P^\alpha}{d\P}\Big|_{\mathcal{F}_t}=Z_\alpha(t).$$
It now follows from Girsanov's theorem that there exist a $\P^\alpha$--standard Brownian motion 
$\{B^\alpha(t),\ t\ge0\}$   such that
\begin{align*}
\int_0^t\sqrt{\frac{M_2(s)}{N}}dB(s)=-\alpha\int_0^t M_2(s)ds+\int_0^t\sqrt{\frac{M_2(s)}{N}}dB^\alpha(s).
\end{align*}
Moreover, we deduce from \eqref{covar} and again Girsanov's theorem (see the statement of Theorem VIII.1.4 p. 327 in \cite{ry}) that for each $k\ge0$ there exists a Brownian motion $\{B^\alpha_k(t),\ t\ge0\}$ with
\begin{align*}
\int_0^t\!\!\sqrt{\frac{X_k(s)(1-X_k(s))}{N}}dB_k(s)\!=\!\int_0^t\!\!\alpha(M_1(s)-k)X_k(s)ds\!+\!\int_0^t\!\!\sqrt{\frac{X_k(s)(1-X_k(s))}{N}}dB^\alpha_k(s).
\end{align*}
Consequently under $\P^\alpha$, we have proved weak existence to our infinite dimensional system
\eqref{un}. We can now turn to the

\noindent{\sc Proof of uniqueness}
We exploit again  Girsanov's theorem to prove weak uniqueness. Consider for some $\alpha, \delta>0$ any $\x_\delta$--valued solution of 
our SDE, which we rewrite as
\begin{align*}
X_k(t)&=x_k+\int_0^t\left[\alpha(M_1(s)-k)X_k(s) +\lambda(X_{k-1}(s)-X_k(s))\right]ds+
\M_k(t),\ k\ge0;\\
M_1(t)&=\sum_{k\ge0}kx_k+\int_0^t\left[\lambda-\alpha M_2(s)\right]ds+\M(t),
\end{align*}
where for $k,\ell\ge0$,
\begin{align*}
\langle \M_k,\M_\ell\rangle_t&=N^{-1}\int_0^tX_k(s)(\delta_{k,\ell}-X_\ell(s))ds,\\
\langle \M_k,\M\rangle_t&=N^{-1}\int_0^tX_k(s)((k-M_1(s))ds,\\
\langle \M,\M\rangle_t&=N^{-1}\int_0^t M_2(s)ds.
\end{align*}
Let $\Q^\alpha$ denote the probability law of our solution
on the space $C([0,+\infty);\x_\delta)$, and define, for $t\ge0$,
$$Y_\alpha(t)=\exp\left(\alpha\sqrt{N}\M_t-\frac{\alpha^2N}{2}\int_0^tM_2(s)ds\right).$$
For each $n\ge1$, let
$$\tau_n:=\inf\left\{t>0,\ \int_0^tM_2(s)ds>n\right\}.$$
It is not hard to show that   for each $n\ge1$, the probability measure $\mathbb{Q}$ defined on $C([0,+\infty);\x_\delta)$ equipped with its Borel $\sigma$--field, by
$$\frac{d\mathbb{Q}}{d\Q^\alpha}\Big|_{\F_{\tau_n}}=Y_\alpha(\tau_n)$$
coincides with the law of the unique weak solution of \eqref{un0} up to time $\tau_n$. Hence 
the restriction of $\Q^\alpha= \left(Y_\alpha(\tau_n)\right)^{-1} \cdot \mathbb{Q}$ to
the $\sigma$--algebra $\F_{\tau_n}$ coincides with the law of the solution which we have constructed above. Since $\tau_n\to\infty$ a.s., weak uniqueness is proved.

\subsection{A comparison theorem for one--dimensional SDEs}
We state a result, which will be useful later in this paper.
 Our processes are defined on a probability space $( \Omega, \mathcal{F},\mathbb{P})$, equipped with a filtration $( \mathcal{F}_t,t\ge 0)$, assumed to satisfy the ``usual hypotheses'', which is such that for each $k, \ell \ge0$ $\left\{B_{k,\ell}(t),\ t\ge0\right\}$ is a $\mathcal{F}_t$--Brownian motion. We denote by $\mathcal{P}$ the corresponding $\sigma$-algebra of predictable subsets of $\R_+ \times \Omega$.

From the weak existence and uniqueness, we deduce that our system has the strong Markov property, using a very similar proof as in Theorem 6.2.2 from \cite{vara}. Indeed, the proof of that results exploits weak uniqueness of the martingale problem, together with the measurability of the law of the solution, with respect to the starting point. In our case that mapping is easily shown to be continuous.

In the next sections, we will use the following comparison theorem several times. This Lemma can be proved exactly as the comparison Theorem 3.7 from chapter IX of \cite{ry}.

\begin{lemma}\label{comparaison}
Let $B(t)$ be a standard $\mathcal{F}_t$--Brownian motion, $T$ a stopping time,  $\sigma$ be a 1/2 H\"older function, $b_1 : \R \to \R$ a Lipschitz function and  $b_2 : \Omega \times \R_+ \times \R \rightarrow \R$ be a $ \mathcal{P}  \otimes B(\R)$ measurable
function. Consider the two SDEs
\begin{equation}\label{Y11}
\left\{
\begin{aligned}
dY_1(t) &= b_1(Y_1(t))dt + \sigma (Y_1(t)) dB(t), \\
Y_1(0)&=y_1;
\end{aligned}
\right.
\end{equation}
\begin{equation}\label{Y12}
\left\{
\begin{aligned}
dY_2(t) &= b_2(t,Y_2(t))dt + \sigma (Y_2(t)) dB(t), \\
Y_2(0)&=y_2.
\end{aligned}
\right.
\end{equation}

Let $Y_1$ (resp $Y_2$) be a solution of (\ref{Y11})  (resp (\ref{Y12})). If $y_1 \le y _2$ (resp $y_2 \le y_1$) and
outside a measurable subset of $\Omega$ of probability zero, $\forall t \in \left[0, T\right]$, $\forall x\in \R,$ $b_1(x) \le b_2(t,x)$ (resp $b_1(x) \ge b_2(t,x)$), then  a.s. $\forall t \in \left[0, T\right],$ $Y_1(t) \le Y_2(t)$ (resp $Y_1(t) \ge Y_2(t)$).
\end{lemma}

\section{The result for a specific set of initial conditions}
From 
$$\P(E) + \P(F) - \P(E \cap F)=\P(E \cup F) \le1,$$ we deduce the following trivial lemma which will be used several times below :

\begin{lemma}\label{trivial}
Let $E,F \in\mathcal{F}$. Then $\P(E \cap F) \ge \P(E) +\P(F)-1$.
\end{lemma}

Now first we show that $M_1$ cannot grow too fast :

\begin{lemma}\label{M_1 pas trop rapide}
For all $c >0$,  $t>0$,  $t'>0$,
$$  \P\left(\sup_{0 \le r \le t' } M_1( t+r) - M_1(t) \le \lambda t' + c \right) \ge 1 - \exp (-2\alpha Nc).$$
\end{lemma}

\bpf
Define $Z^t_{t+s}= \int_t^{t+s} \sqrt{\frac{M_2(r)}{N}} dB_r - \alpha \int_t^{s+t} M_2(r) dr $.
We note that, for any $t>0$, $\{\exp( 2 \alpha N Z^t_{t+u}),\  u \ge 0\}$ is both a local martingale and a super--martingale. We also have 
$$\sup_{0 \le s \le t'} M_1( t+s) - M_1(t)  \leq \sup_{0 \le s \le t' } Z^t_{t+s} + \lambda t'.$$
But for all $c >0$,
\begin{align*}
 \P\left( \sup_{0 \le u \le t' } Z^t_{t+u} \ge c\right) &\le  \P\left( \sup_{0 \le u \le t'  } \exp\left( 2 \alpha N Z^t_{t+u}\right) \ge \exp\left( 2 \alpha N c\right)\right) \\
 &\le \exp \left(-2\alpha Nc\right), 
 \end{align*}
 where we have taken advantage of the fact that $\exp\left( 2\alpha N Z^t_{t+u}\right) $ is a local martingale and of Doob's inequality. 
Then
$$\P\left(\sup_{0 \le r \le t' } M_1( t+r) - M_1(t) \le \lambda t' + c \right) \ge 1 - \exp \left(-2 \alpha Nc\right).$$
\epf

Note that we have in fact 
$ \P\left( \sup_{u \ge 0 } Z^t_{t+u} \ge c \right) \le \exp \left(-2 \alpha Nc\right)$.

We choose an arbitrary value $m>0$ for $M_1(0)$, which will remain the same throughout this document (for example one could choose $m=1$), and we define 
\begin{equation}\label{eq1}
{\overline{\ep}}=\frac{1}{10N\alpha},\quad{t'_3} = \frac{{\overline{\ep}} N} {3 \lambda }= \frac{1}{30 \lambda \alpha},
\end{equation}
$$m_{\max} = m + \lambda A({t'_3}) + \frac{{\overline{\ep}}}{6},$$ where $A(t)= \frac{1}{4N} \int_0^t (1-X_0(s) ) ds$, 
\begin{equation}\label{eq2}
 p_2 =   \exp(-\alpha N\frac{{\overline{\ep}}}{6})= \exp(-\frac{1}{60}), 
\end{equation} 
\begin{equation}\label{eq3}
\mu =\frac{ {\overline{\ep}}} {6 m_{\max}}\wedge \frac{{\overline{\ep}}}{4} \wedge \frac{1}{10},
\end{equation}
 and let $\delta$ be a real number, which will be specified below, such that $\delta \le \frac{1}{10} \wedge \frac{{\overline{\ep}}}{m}$.

Now let ${Y_0}$ be the solution of the following SDE :
\begin{equation}\label{eqref}
\left\{
\begin{aligned}
d{Y_0}(t)&=  dt + 2 \sqrt{{Y_0}(t)} dW(t)\\
{Y_0}(0)&= \delta
\end{aligned}
\right.
\end{equation}
with $W$ a standard Brownian motion.

We will show that starting with $X_0 (0)=x_0 \le \delta$ , $ M_1(0)=m_1 \le m$, and as long as $X_0 M_1 < 2{\overline{\ep}}$ and $X_0$ remains small enough, we can compare $X_0(t)$ with the solution of (\ref{eqref}). 

\begin{lemma} \label{simple model}
For any $\delta>0$, ${\overline{\ep}}$ given as in \eqref{eq1}, $\mu$ as in \eqref{eq3},
let $$T_{min} = \inf \{t>0, X_0(t) M_1(t) > 2 {\overline{\ep}} \text{ or }X_0(t) > \delta + \mu \} .$$ Then provided that $X_0(0)=x_0 
\le \delta$, if  $A(t):= \frac{1}{4} \int_0^t \frac{1-X_0(s)}{N}ds$, 
there exists a standard Brownian motion $W$ such that the corresponding solution $Y_0$ of 
\eqref{eqref} satisfies  
$$X_0(t) \le {Y_0}(A(t)),\quad \forall t \in \left[ 0,T_{min} \right] .$$ 
\end{lemma}

\bpf
 We first note that for $0\le t \le T_{min}$,  $\frac{t}{5N} \leq A(t) \leq \frac{t}{4N}$ because $\frac{4}{5} \leq 1- X_0(t) \leq 1 $ (thanks to the choices of $\mu$ and $\delta$, and $1- X_0(t) \geq 1- \delta - \mu \geq  1 - \frac{1}{10} - \frac{1}{10} \geq \frac{4}{5}$). 
 
Define $\sigma(t)=\inf \{ u>0, A(u) \ge t \}$ and $\tilde{X}_0(t)=X_0(\sigma(t))$ ( resp $\tilde{M}_1(t)=M_1(\sigma(t))$). Then there exists a standard Brownian motion $W_t$ such that
 
\begin{align*}
d\tilde{X_0} ( t) =(  \alpha \tilde{M_1}(t) - \lambda  ) \tilde{X_0} (t) \frac{4N}{1-\tilde{X_0}(t)}dt +2 \sqrt{\tilde{X_0}(t)}dW_t\, .
\end{align*}
But whenever $t\le A(T_{\min})$,
\begin{align*}
( \alpha \tilde{M_1}(t) - \lambda  )\tilde{X_0} (t) \frac{4N}{1-\tilde{X_0}(t)} &\le
 \frac{4\alpha N\tilde{M_1}(t)\tilde{X_0} (t)}{1-\tilde{X_0}(t)}\\
 &\le 1,
 \end{align*}
 because the numerator on the right is less than or equal to $4/5$, while the denominator is bigger than or equal to the same figure.

The result then follows from Lemma \ref{comparaison}.
\epf

Next we will prove that ${Y_0}$ reaches zero with positive probability on a fixed time interval. 
For any $\alpha\in\R$, we define
$$T'_\alpha=\inf\{t>0,\ Y_0(t)=\alpha\}.$$
\begin{lemma}\label{equation reference}
Let $\{Y_0(t),\,t\ge0\}$ be the solution of (\ref{eqref}). For all $p <1 $, $\tilde{\mu} >0$,  there exists $\delta >0$ such that with $t'_3$ defined as in \eqref{eq1},
$$\P(T'_0 \le {t'_3}   \wedge T'_{\delta+\tilde{\mu}})  \ge p.$$
\end{lemma}

\bpf
Let
\begin{align*}
\tilde{Y}(t)&=\delta \exp\left(  -t + 2 W (t)\right),\\
 D(t)&=\int_0^t \tilde{Y}(s) ds, \\
\rho (t) &= \inf \{s>0, D(s) >t\}.
\end{align*}
It is not too hard to show that there exists a Brownian motion $\bar W$ such that
$$\tilde{Y}(\rho(t))=\delta\exp\left(-\int_0^t\frac{ds}{\tilde{Y}(\rho(s))} +2\int_0^t
\frac{d\bar W(s)}{\sqrt{\tilde{Y}(\rho(s))}}\right).$$
It now follows from Ito's formula that the process $\{Y(t):=\tilde{Y}(\rho(t)),\ t\ge0\}$ is the unique strong solution of equation \eqref{eqref} driven by $\bar W$, hence $Y_0(t)=\tilde{Y}(\rho(t))$, 
$t\ge0$. We deduce that $T'_0=D(\infty ) < \infty$, and
\begin{align*}
&\P(T'_0 \le {t'_3}   \wedge T'_{\delta+\tilde{\mu}})   \\
&=  \P \left(  \left\{ \int_0^\infty  \exp(  -t + 2 W (t) )dt \le \frac{ {t'_3}}{\delta} \right\} \cap \left\{ \sup_{t\ge 0} \exp( - t +2W(t)) \le \frac{\delta + \tilde{\mu}}{\delta} \right\} \right) \\
&\to 1,
\end{align*}
as $\delta \rightarrow 0$, since $\sup_{t\ge 0} \exp( - t +2W(t)) < \infty$ a.s.

\epf

Now we can choose the value of $\delta$ which we will be using from now on. Let $\delta'$ be the largest value of $\delta$ such that Lemma \ref{equation reference} holds, with   $p=p_2$ defined in
\eqref{eq2} and $\tilde{\mu}=\mu$ ( which is a function of $m_{\max}$) as defined by \eqref{eq3}. We choose (recall that the value of ${\overline{\ep}}$ has been defined in \eqref{eq1})
 \begin{equation}\label{choice-delta}
 \delta =  \delta' \wedge \frac{1}{10} \wedge \frac{{\overline{\ep}}}{m}.
 \end{equation}

Thanks to Lemma \ref{equation reference}, when starting at time 0 from $\delta$, ${Y_0}$  reaches 0 with probability $p_2$ before time ${t'_3} \wedge T'_{\delta+\mu}$. Then $X_0$ will do the same before time $A({t'_3}) \wedge A(T'_{\delta+\mu})$, provided that $X_0(t) M_1(t) \leq 2 {\overline{\ep}} $,  $\forall 0\le t \le A({t'_3}) \wedge A(T'_{\delta+\mu})$. Hence the fact that $T_0 < A({t'_3})$ with positive probability, provided $x_0 \le \delta$ and $M_1(0) \le m$  will follow from the above results and

\begin{lemma} \label{M_1 step 2}
If $X_0(0) \le \delta$ and $M_1(0) \le m$, then we have (again with ${\overline{\ep}}$ and $t'_3$
given by \eqref{eq1})
$$\P\left(\sup_{ 0 \le t \le A({t'_3} \wedge T'_{\delta+\mu})} X_0 (t) M_1 (t) \le 2 {\overline{\ep}}\right)=p_3 > 1- p_2.$$
\end{lemma}
\bpf
We use Lemma \ref{M_1 pas trop rapide}.
Consider the event $$E_{m,{t'_3}, \bar{\ep}}=\left\{ \sup_{ 0 \le t \le A({t'_3}) \wedge A(T'_{\delta+\mu})} M_1(t)  \le m + \lambda A({t'_3}) + \frac{{\overline{\ep}}}{6} \right\}.$$
We have
\begin{align*}
\P( E_{m,{t'_3}, \bar{\ep}}) &\geq  \P\left(\sup_{ 0 \le t \le A({t'_3})}  M_1 (t) \le m_1 + \lambda A({t'_3}) + \frac{{\overline{\ep}}}{6} \right)\\
&\ge 1 - \exp(- \alpha N \frac{{\overline{\ep}}}{3})= 1- \exp(-\frac{1}{30}).
\end{align*}
Since $X_0(t) \le \delta + \mu$ for $t \le A(T'_{\delta+\mu})$,  on the event $E_{m,{t'_3}, \bar{\ep}}$, 

\begin{align*}
 \sup_{ 0 \le t \le A({t'_3}) \wedge A(T'_{\delta+\mu})} X_0 (t) M_1 (t)  & \le (\delta + \mu) ( m + \lambda A({t'_3}) + \frac{{\overline{\ep}}}{6} ) \\
& \le \delta m + \mu m+ \lambda A({t'_3})+  \frac{{\overline{\ep}}}{6} \\
& \le {\overline{\ep}} + \frac{{\overline{\ep}}}{6} + \frac{{\overline{\ep}}}{12} + \frac{{\overline{\ep}}}{6} \\
& \le 2 {\overline{\ep}},
\end{align*}
where we have used the fact that $\delta + \mu \le 1$ for the second inequality.
\epf


Combining Lemma \ref{trivial}, Lemma \ref{simple model}, Lemma \ref{M_1 step 2} and Lemma \ref{equation reference}, 
denoting $t_3=A(t'_3)$, we deduce the

\begin{corollary}\label{departpetitset}
There exists $p_{\text{fin}}\ge p_3 + p_2 -1>0$ such that 
$$\P(T_0 \le  t_3 | X_0(0) \le \delta, M_1(0) \le m )  \ge p_{fin} >0.$$
\end{corollary}

While this Corollary is rather intuitive, we shall need the slightly more general following result, i.e. 
with a larger set of initial conditions. Given $m$ as above, and $\delta$  as in \eqref{choice-delta}, let 
$$\mathcal{I}= \{ (x_0,m_1) \in [ 0, 1] \times \R_+,  x_0 \le \delta, x_0 m_1 \le \delta m \}.$$
   We now prove the ($p_{\text{fin}}$ is as defined in Corollary \ref{departpetitset})
   \begin{proposition}
   For any $(x_0,m_1)\in\mathcal{I}$, 
   $$\P(T_0\le t_3|X_0(0)=x_0,M_1(0)=m_1)\ge p_{\text{fin}}.$$
   \end{proposition}
   \bpf
 Thanks to the previous Corollary, we only need to consider the case  $m_1 > m$.
Let $(x_0, m_1)$ be a point in the set $\mathcal{I}$. First, let us consider the point $(\delta, m)$. From the previous section, starting from $(\delta,m)$, the process $(X_0,M_1)$ has a strictly positive probability to reach $0$ before time ${t_3}=A({t'_3})$. We will show that the process starting from $(x_0,m_1)$ has a larger probability to reach $0$ before time ${t_3}$, which will prove the Proposition.

Let $C= \frac{m_1}{m} \geq 1$. Then we have $ x_0 \leq \frac{\delta}{C} $.

Now we will use the same reasoning as in Lemma \ref{equation reference} with a few modifications. Indeed, since the probability that $ {Y_0}(t)$ reaches 0 before a prescribed time is decreasing in $\delta$, we increase this probability by starting from ${Y_0}(0)=x_0= \delta'\leq \frac{\delta}{C} $, since $C \geq 1$.  We will use this new value. Moreover, the starting point satisfies $x_0 m_1 \leq  {\overline{\ep}}$. The only thing which is worse than with the starting point $(\delta,m)$ is the fact that $m_1$ is greater than $m$, hence a greater $m_{\max}$. But this only appears in one place : in the definition of $\mu$.

Note that if we define $m'_{\max}= m_1 + \lambda {t_3} + \frac{{\overline{\ep}}}{6}$,  we deduce from
Lemma \ref{M_1 pas trop rapide}
\begin{align*}
\P( \sup_{0 \le t \le {t_3}} M_1(t)  \le m'_{\max}  )\geq 1 - \exp(- \alpha N \frac{{\overline{\ep}}}{3}).
\end{align*}

We define $\mu'$ similarly as $\mu$ in \eqref{eq3}, but with $m_{\max}$ replaced by $m'_{\max}$, hence since $m'_{\max}\le C m_{\max}$, $\mu' \ge \frac{\mu}{C}$.  But if we look at the proof of Lemma \ref{equation reference}, we have, since $\frac{ {t'_3}}{\delta'} \ge \frac{ C{t'_3}}{\delta} \ge \frac{ {t'_3}}{\delta}$ and $\frac{\delta' + \mu'}{\delta'} = 1 + \frac{\mu'}{\delta'} \ge 1 + \frac{\mu}{\delta}$,

\begin{align*}
\P(T_0\le t_3)&\ge
 \P(T'_0 \le {t'_3}   \wedge T'_{\delta'+\mu'})\\  &\ge  \P \left(  \left\{ \int_0^\infty  \exp(  -t + 2 W (t) )dt \le \frac{ {t'_3}}{\delta'} \right\} \cap \left\{ \sup_{t\ge 0} \exp( - t +2W(t)) \le \frac{\delta' + \mu'}{\delta'} \right\} \right) \\
 & \ge  \P \left(  \left\{ \int_0^\infty  \exp(  -t + 2 W (t) )dt \le \frac{ {t'_3}}{\delta} \right\} \cap \left\{ \sup_{t\ge 0} \exp( - t +2W(t)) \le \frac{\delta+ \mu}{\delta} \right\} \right)
\end{align*}

Hence we have a larger probability to reach zero starting from $(x_0, m_1)$ rather than from $(\delta,m)$, which concludes the proof.

\epf

We sum up in the following Proposition the results obtained in this section, with $\ep= \delta m$ (recall that $m$ has been chosen arbitrarily, $\delta$ is prescribed by \eqref{choice-delta}, and note that $\ep \le \overline{\ep}$). 

\begin{proposition}\label{clef}
Let $X(t)= (X_k(t))_{k \in {\Z_+}}$ be the solution of \eqref{un} , and $M_1(t)=\sum_{k\ge1}kX_k(t)$. Then there exist
 $p_{fin} >0$ and ${t_3}$ such that for any $t\ge0$,
$$\P(T_0\le t+{t_3} | X_0(t) \le \delta,X_0(t) M_1(t) \le \ep ) \ge p_{fin}>0.$$
\end{proposition}

\section{A recurrence property of $M_1$}
With the help of the results proved in the previous section, we will now prove some results on
 $M_1$. We will show that as long that as the ratchet has not clicked, $M_1$ is bound to return under some specified value. This particular point will be  important in the sequel.

We begin with the following lemma, which is true for any probability on ${\Z_+}$. It will be crucial for establishing one of our first estimates.

\begin{lemma}\label{minorationUn}
Let p be a probability on ${\Z_+}$, and let $x_k=p(k)$, $m_1 = \sum_{k \ge0} k x_k$ and $m_2 = \sum_{k \ge0} (k- m_1)^2 x_k $. Then 
$$m_2\ge (1-x_0 )m_2 \ge  x_0 m_1 ^2.$$
\end{lemma}

\bpf
If $x_0=1$, $m_1=m_2=0$ and the result is true. So it suffices to study the case $x_0<1$. By Jensen's inequality we have
\begin{align*}
 \left(\sum_{k \ge 1 } \frac{x_k}{1-x_0}k\right)^2 \le \sum_{k \ge 1 } \frac{x_k}{1-x_0}k^2
\end{align*}
with equality if and only if there exists only one $k \ge 1$ such that $x_k > 0$. Then :

\[
 \left(\sum_{k \ge 1 } x_k k\right)^2 \le (1-x_0)\sum_{k \ge 1 } x_k k^2, 
 \]
that is\[
 m_1^2 \le (1-x_0)\sum_{k \ge 1 } x_k k^2 ,
 \]
 hence
 \begin{align*}
 x_0 m_1^2 &\le (1-x_0)\sum_{k \ge 1 } k^2 x_k   - (1-x_0)m_1^2 \\
&= (1-x_0) m_2 .
\end{align*}~\epf

We now introduce for each $t\ge0$ the stopping time
$${H^t_\lambda}:=\inf\{s \ge t,\ X_0(s )M_1(s)^2\le 2\frac{\lambda+1}{\alpha}\},$$
and we define $H_\lambda = H_\lambda ^0$.

Our next claim is 
 \begin{proposition}\label{rec-XM2} 
 For any stopping time $T$, we have $H^T_\lambda < + \infty$ a.s.
\end{proposition}

The Proposition follows from the strong Markov property and
\begin{lemma}\label{return}
Suppose that $X_0(0)M_1(0)^2>2\frac{\lambda+1}{\alpha}$. Then ${H_\lambda}<\infty$ a.s.
\end{lemma}
\bpf
On the interval $[0,{H_\lambda}]$,we have from Lemma \ref{minorationUn} $$ - \frac{\alpha}{2}M_2\le - \frac{\alpha}{2} X_0 M_1^2 \le -( \lambda +1),$$  hence from the third line of \eqref{3eq}, 
\begin{equation}\label{sdeY}
\begin{aligned}
M_1(t)  
 \le  M_1(0) -t -\frac{\alpha}{2}\int_0^t M_2(s) ds+ \int_0^t \sqrt{\frac{M_2(s)}{N}}dB_s.
\end{aligned}
\end{equation}
We will show next that $$Z_t:=\int_0^t\sqrt{\frac{M_2(r)}{N}}dB_r-\frac{\alpha}{2}\int_0^tM_2(r)dr$$ is bounded from above a.s. This will imply the result, since on the event $\{H_\lambda=+\infty\}$,
\eqref{sdeY} holds for all $t>0$, which would imply that $M_1$ eventually becomes negative, and this is absurd.

If we define $ C(t)=\frac{1}{N} \int_0^tM_2(s)ds$, we have $Z_t=W(C(t))-\frac{\alpha N}{2} C(t)$ where $W$ is a standard Brownian motion. 

Now,  if $C(\infty)=\infty$ then $\lim_{t\to\infty} Z_t=-\infty$, hence $Z_t$ is bounded from above. Or else $ C(\infty)<\infty$, and we have  
$\sup_{t>0} |  Z_t |  =\sup_{0<s<C(\infty)} | W(s)-\frac{\alpha N}{2} s | < \infty$ a.s.
\epf

Now we will finally be able to prove that $M_1$ always returns below $\beta:=\lambda/\alpha$, as long as the ratchet does not click.
Let for each $t\ge0$
$$S^t_\beta = \inf\{s >t , M_1(s) \leq \beta \}.$$ 

We have the following lemma :

\begin{lemma}\label{Mini} For any stopping time $T$,
 $\P(T_0 \wedge S^T_\beta < \infty)=1$.
\end{lemma}

\bpf From the strong Markov property of the solution of
\eqref{un}, we may assume that $T=0$. First, we let $\delta_{\inf}=\delta \wedge \frac{\ep^2 \alpha}{4 (\lambda +1)}$ (recall that $\ep=\delta m$).

Now we introduce the process $Y^{s}_{{t}}$, defined for all $s \ge 0 \text{, }  {{t}} \geq s$, which is the solution of the following SDE :

\begin{equation}
\left\{
\begin{aligned}
 dY^{s}_{{t}} &=  \sqrt{\frac{Y^{s}_{{t}}(1-Y^{s}_{{t}})}{N}}dB_0({t}),  {t} \ge s, \\
 Y^{s}_{s}&= \delta_{\inf}.
\end{aligned}
\right.
\end{equation}

We define for any $0 \le u \le 1$ $$R^{s}_u=\inf\{{{t}} \geq s, Y^{s}_{{t}}=u \}.$$
We have 
\begin{equation}\label{Yexit}
\left\{
\begin{aligned}
&R^{s}_0 \wedge R^{s}_1 < + \infty \text{ a.s.},\\
&\P(R^{s}_1<R^{s}_0)>0.  
\end{aligned}
\right.
\end{equation}
 Indeed, for all $a \in \left( 0,\delta_{\inf}\right)$, by the non--degeneracy of the diffusion coefficient, $Y^s_{{t}}$ gets out of $\left[ a, 1-a \right]$ in finite time. Then if we choose $a$ small enough (using the same reasoning as in Lemma \ref{equation reference}), we have a chance $p'_{fin}$ to reach $0$ before a time $V>0$ as soon as we start below $ a$ (the same with $1$ and starting above $\ge 1-a$ by symmetry).

Define recursively the stopping times
\begin{align*}
\xi_1&=\inf\{t>s;\, Y^s_t\not\in(a,1-a)\},\\
\text{and for }k\ge1,\ \xi_{k+1}&=\inf\{t>\xi_k+V;\, Y^s_t\not\in(a,1-a)\}.
\end{align*}
A standard application of the strong Markov property of $Y^s_t$ yields that
$$\P(Y^s_{\xi_{k+1}}\not\in\{0,1\})\le(1-p'_{fin})^k,$$
hence the first line of \eqref{Yexit}. The second line is essentially obvious.
Note that using an argument based upon Green's functions, one can in fact prove that $\E(R^{s}_0 \wedge R^{s}_1)  <+ \infty$.

From this we deduce that there exist $K>0$, $p>0$ such that $ \P(R^{s}_1 \leq K \wedge R^{s}_0) \geq p >0$. In particular $\P(R^{s}_1 \leq K ) \geq p >0$.
 
 We define $L=K \vee {t_3}$, where ${t_3}$ is as in Proposition \ref{clef}, and the following collection  of stopping times : 
 \[U^t_0=\inf \left\{s>t,  X_0(s) M_1^2(s) \le 2 \frac {\lambda +1}{\alpha}\right\},\] 
 and for all $n\geq 1$,  
 \[
 U^t_n=\inf \left\{s>U^t_{n-1}+L,  X_0(s) M_1^2(s) \le 2 \frac {\lambda +1}{\alpha}\right\}.
 \] 
For all $n\ge0$, $U^t_n$ is a.s. finite, thanks to Proposition \ref{rec-XM2}.
 
 Now, at time $U^t_0$ :
either we are on the event $A_0=\{X_0(U^t_0) \leq \delta_{\inf}\}$   $( \delta_{\inf}\leq \delta)$, in which case
\begin{align*}
X_0(U^t_0) M_1(U^t_0) &= \sqrt{X_0(U^t_0)M_1^2(U^t_0) \times X_0(U^t_0)}\\
&\le \sqrt{2\frac{\lambda+1}{\alpha} \frac{\ep^2 \alpha}{4(\lambda+1)}}\\
& < \ep,
\end{align*}
and we deduce from Proposition \ref{clef} that 
$$\P(T_0 \le U^t_0+L|A_0)=p_{fin} >0.$$

The other possibility is that we are on the event 
\begin{align*}
A_0^c&=B_0\cup C_0,\qquad\text{where}\\
B_0&=\{X_0(U^t_0) > \delta_{\inf}\}\cap\{\inf_{ U^t_0 \leq s \le U^t_0+L} M_1(s) \geq \beta\},\\
C_0&=\{X_0(U^t_0) > \delta_{\inf}\}\cap\{\inf_{ U^t_0 \leq s \le U^t_0+L} M_1(s) < \beta\}.
\end{align*}
On the event $C_0$, $S_\beta\le U^t_0+L$. 
On the event $B_0$, 
$\inf_{U^t_0\le s\le U^t_0+L}(\alpha M_1(s) - \lambda )X_0(s) \geq 0$, and then we deduce from Lemma \ref{comparaison}
 that $X_0(s) \ge Y^{U^t_0}_s$. Consequently, if  $T_1=\inf\left\{ t \ge 0, X_0(t)=1\right\}$,
 $$\P(T_1 \le U^t_0+L|B_0) \ge p >0 .$$ But if $X_0(s)=1$, then $M_1(s)=0$. Hence 
 $$\P(S_{\beta} \le U_0^t+L|B_0) \ge p >0.$$ Finally
 \begin{align*}
 \P(T_0 \wedge S^t_\beta \le U^t_0+ L)&=
 \P(T_0 \wedge S^t_\beta \le U^t_0+ L|A_0)\P(A_0)+\\
 &\quad \P(T_0 \wedge S^t_\beta \le U^t_0+ L|B_0)\P(B_0)+
   \P(T_0 \wedge S^t_\beta \le U^t_0+ L|C_0)\P(C_0)\\
   &\ge p_{fin}\P(A_0)+p\P(B_0)+\P(C_0)\\
   &\ge p_{fin}\wedge p=:q,
   \end{align*}
 $$\P(T_0 \wedge S^t_\beta = + \infty) \le \P(T_0 \wedge S^t_\beta \ge U^t_0+ L) \le 1-q.$$

It follows from the strong Markov property of the process $ X=(X_k,\ k\ge0)$, repeating this argument with $U^t_0$ replaced by $U^t_1$ that
$$\P(T_0 \wedge S^t_\beta = + \infty) \le \P(T_0 \wedge S^t_\beta \ge U^t_1+ L) \le (1-q)^2.$$

Iterating the above argument, we have for all $\ell\ge0$, 
$$\P(T_0 \wedge S^t_\beta > U^t_\ell +L) \le (1-q)^{\ell}.$$ 
We have proved that
 $$\P(T_0 \wedge S^t_\beta = + \infty) =0.$$
\epf

\section{Reaching the special set from any initial condition}

Now we will show that starting from an initial condition $((x_k)_{k \in {\Z_+}},m_1)$ with $ m_1 \leq \beta$ the process has a   probability bounded below by $p_{fin}$ to click before a given time. Since the process is Markovian and this situation repeats itself as long as the ratchet has not clicked, we will conclude that $\P(T_0< + \infty) =1$.

In this section we denote by $(x_k)_{k \ge0}$ the initial condition of our system, and we suppose that $m_1=\sum_{k\ge0}kx_k \leq \beta$.

One of the difficulties we have to face  is that the quadratic variation of $X_0$ is $\frac {X_0(1-X_0)}{N}$, which is not bounded away from $0$, near $1$ and $0$. We need to study three separate cases. 

The first case will be described in terms of the constant
\begin{equation}\label{x_max}
{x_{\max}}= \max \left\{\frac{9}{10} ,  \frac{3 \lambda+ 5\alpha }{5 (\lambda + \alpha) }, 1- \frac{2}{\lambda} \right\}. 
\end{equation}

\subsection{$x_0 \in \left( {x_{\max}}; 1\right]$}

The following lemma will show that if $X_0$ starts close to $1$, it will quickly go under ${x_{\max}}$ :

\begin{lemma}\label{X01}
Let $t_1=\frac{8}{ \lambda^2}$. If $X_0(0) >{x_{\max}}$ , then 
$$ \P (\inf_{s<t_1} X_0 (s) \le {x_{\max}}) \geq 1- \exp(-N).$$

\end{lemma}

\bpf
Let $T_{{x_{\max}}}=\inf \{s\ge 0, X_0 (s) \le {x_{\max}}  \}  $. On the time interval $[0,T_{{x_{\max}}})$, we have$$X_0(s) >{x_{\max}} \ge  \frac{3 \lambda+ 5\alpha }{5 (\lambda + \alpha) }. $$ 
Since $X_1 \le 1 -X_0$, on the same interval we have $X_1(s) \le \frac{2 \lambda}{5( \lambda + \alpha)} $ while $X_0(s)>\frac{9}{10}$, hence
\begin{align*}
  \alpha M_1(s) X_1(s) +\lambda X_0(s) - (\lambda +\alpha) X_1(s)    \geq &    
  \lambda X_0(s) - (\lambda +\alpha) \frac{2 \lambda }{5 (\lambda + \alpha) }\\
  \geq &  \frac{\lambda}{2}.
\end{align*}
Then $X_1(s) \geq Y_1(s)$ for $s \in \left[0,T_{{x_{\max}}}\right]$, where $Y_1$ is the solution of the SDE
\begin{equation}
\left\{
\begin{aligned}
dY_1(s) &=\frac{\lambda} {2}ds + \sqrt{\frac{Y_1(1-Y_1)}{N}} dB_1(s),\\
Y_1(0)&=0,
\end{aligned}
\right.
\end{equation}
where we stop $Y_1$ as soon as it reaches $1$.

 We have 
\begin{align*}
&\P \left( \int_{0}^{t_1}\sqrt{\frac{Y_1(1-Y_1)}{N}} dB_1 < -C \right) \\
&=\P \left( - \int_{0}^{t_1}\sqrt{\frac{Y_1(1-Y_1)}{N}} dB_1 > C \right)  \\
 & \leq \P \left(  \exp \left( -\gamma\int_{0}^{t_1}\sqrt{\frac{Y_1(1-Y_1)}{N}} dB_1 - \int_0^{t_1} \frac{\gamma^2 Y_1(s) (1-Y_1(s))}{2N} ds \right)  >  \exp \left( \gamma C - \frac{\gamma^2}{8N} t_1 \right) \right)\\
& \leq \exp \left( -\gamma C + \frac{\gamma^2}{8N} t_1 \right) ,
\end{align*}
where the first inequality follows from $Y_1(s)(1-Y_1(s))\le1/4$, and the second one is Chebychev's inequality.
Choosing $\gamma=4CN/t_1$ and $C=2/\lambda$, we deduce that
$$\P\left(\int_{0}^{t_1}\sqrt{\frac{Y_1(1-Y_1)}{N}} dB_1 \geq -\frac{2}{\lambda}  \right) \geq 1 - \exp \left(-  N   \right) >0.$$

Now, since 
$$\int_{0}^{t_1} \frac{\lambda} {2} ds= \frac{4}{ \lambda},$$
and on $[0,T_{x_{\max}})$, $X_0(s)>1-2/\lambda$, hence $X_1(s)<2/\lambda$, we
have the inclusion
$$  \left\{\int_{0}^{t_1}\sqrt{\frac{Y_1(1-Y_1)}{N}} dB_1 \geq -³\frac{2}{\lambda}\right\}  \subset \{ T_{{x_{\max}}} < t_1 \},$$
which implies that 
$$\P \left( T_{{x_{\max}}} \le t_1 \right) \ge 1- \exp(-N), $$
hence the conclusion.
\epf 

 We need to control $M_1$ on the same time interval of length $t_1$. Using Lemma \ref{M_1 pas trop rapide} we will deduce the following Proposition :

\begin{proposition}\label{stepun}
Let again ${x_{\max}}$ be given by \eqref{x_max}, $t_1=\frac{8}{ \lambda^2}$, $\ep_0 = \frac{1}{2 \alpha N } \ln \left( \frac{2}{1-\exp(-N)} \right)$ and $\beta'= \beta + \lambda t_1+ \ep_0$. If $X_0(0) >{x_{\max}}$ and $M_1(0) < \beta$ , then
$$ \P \left( \{ T_{{x_{\max}}} \le t_1\} \cap \{M_1(T_{{x_{\max}}}) \le \beta'\} \right)  = p_{init}  >0.$$
\end{proposition}

\bpf
It follows from Lemma \ref{X01} and Lemma \ref{M_1 pas trop rapide} 
\begin{align*}
\P(T_{{x_{\max}}} \le t_1) &\ge 1- \exp(-N),\\
\P( M_1(T_{{x_{\max}}}) &\le \beta')\ge 1- \exp(-2\alpha N\ep_0).
\end{align*}
 Those two inequalities together with Lemma \ref{trivial} imply
\begin{align*}
 \P \left( \{ T_{{x_{\max}}} \le t_1\} \cap \{M_1(T_{{x_{\max}}}) \le \beta'\} \right)
& \ge 1-\exp(-N) -\exp( - 2 \alpha N \ep_0)  \\
& = \frac{1- \exp(-N)}{2} =: p_{init}.
\end{align*}
\epf

So even if we started with $(X_0(0),M_1(0))$ such that $X_0(0) > {x_{\max}}$ and $M_1(0) < \beta$, we obtain before time $t_1$  with 
 probability at least
 $p_{init}>0$ a new initial condition $X_0 \leq {x_{\max}}$ and $M_1 \leq \beta'$, so we can resume with the next case.

\subsection{$X_0 \leq  {x_{\max}} $ but either $X_0 > \delta$ or $X_0 M_1 > \ep$}

The idea of this subsection is to show that with a strictly positive probability $p_{trans} $, both
$X_0$ goes from ${x_{\max}}$ to  a $\delta' < \delta$ in finite time, and during the same time interval, $M_1$ stays small enough so that at the end $X_0M_1 \leq \ep$.

We start by showing some inequalities.

\begin{lemma}\label{Bdescend}
Let $\{V_t,\ t\ge0\}$ be a standard Brownian motion, and $c>0$ a constant. 
Then  for any $t>0$, $\tilde{\delta}>0$, $\tilde{\mu}>0$,
\begin{align*}
\P&\left(\inf_{0\le s\le t}\{cs+V_s\}\le-\tilde{\delta}, \sup_{0\le s\le t}\{cs+V_s\}\le \tilde{\mu}\right)\\
&\quad\ge1-\sqrt{\frac{2}{\pi}}\left(\frac{\tilde{\delta}}{\sqrt{t}}+c\sqrt{t}\right)-
2\exp\left[-\frac{1}{2}\left(\frac{\tilde{\mu}}{\sqrt{t}}-c\sqrt{t}\right)^2\right].
\end{align*}
\end{lemma}
\bpf 
Using Lemma \ref{trivial}, the result follows from the two following computations. We have, with $Z$ denoting a $N(0,1)$ random variable,
\begin{align*}
\P\left(\inf_{0\le s\le t}\{cs+V_s\}\le-\tilde{\delta}\right)&\ge\P\left(\inf_{0\le s\le t} V_s\le-\tilde{\delta}-ct\right)\\
&= \P\left(\sup_{0\le s\le t} V_s\ge\tilde{\delta}+ct\right)\\
&= 2\P(V_t\ge\tilde{\delta}+ct)\\
&= 1-\P\left(|Z|\le \frac{\tilde{\delta}}{\sqrt{t}}+c\sqrt{t}\right)\\
&\ge1-\sqrt{\frac{2}{\pi}}\left(\frac{\tilde{\delta}}{\sqrt{t}}+c\sqrt{t}\right).
\end{align*}

On the other hand,
\begin{align*}
\P\left(\sup_{0\le s\le t}(cs+V_s)\le\tilde{\mu}\right)&\ge\P\left(\sup_{0\le s\le t}V_s\le\tilde{\mu}-ct\right)\\
&=1-\P\left(\sup_{0\le s\le t}V_s\ge\tilde{\mu}-ct\right)\\
&=1-2\P\left(Z\ge\frac{\tilde{\mu}}{\sqrt{t}}-c\sqrt{t}\right),
\end{align*}
and
\begin{align*}
\P\left(Z\ge\frac{\tilde{\mu}}{\sqrt{t}}-c\sqrt{t}\right)
&=\P\left(\exp(\gamma Z-\gamma^2/2)\ge\exp\left(\gamma\left[\frac{\tilde{\mu}}{\sqrt{t}}-c\sqrt{t}\right]
-\frac{\gamma^2}{2}\right)\right)\\
&\le\exp\left(-\gamma\left[\frac{\tilde{\mu}}{\sqrt{t}}-c\sqrt{t}\right]+\frac{\gamma^2}{2}\right).
\end{align*}  
Choosing $\gamma=\tilde{\mu}/\sqrt{t}-c\sqrt{t}$, we conclude from the above computations that
$$\P\left(\sup_{0\le s\le t}(cs+V_s)\le\tilde{\mu}\right)\ge1-2\exp\left[-\frac{1}{2}\left(\frac{\tilde{\mu}}{\sqrt{t}}-c\sqrt{t}\right)^2\right].$$
\epf

We will choose from now on
\begin{equation}\label{ep}
\tilde{\ep}=\frac{\log(4)}{2\alpha N},\quad\text{so that }e^{-2N \alpha\tilde{\ep}}=\frac{1}{4}.
 \end{equation}

 We start from $(X_0,M_1)=(x,\beta')$, where $0<x\le{x_{\max}}<1$ (recall the definition \eqref{x_max} of $x_{\max}$) and $\beta<\beta'$. 
Let  $0<\tilde{\mu}=\frac{1-{x_{\max}}}{2}$. We are going to prove that, with positive probability,
$X_0$ goes down to $\delta'$ in a finite number
of steps, while staying below $x+\tilde{\mu}$ (so that $1-X_0(t)\ge a:=\frac{1-{x_{\max}}}{2}$), and while $M_1$ remains under control. 

Considering the SDE
$$dX_0(t)=(\alpha M_1(t)-\lambda)X_0(t)dt+\sqrt{\frac{X_0(t)[1-X_0(t)]}{N}}dB_0,$$
let 
\begin{align*}
A(t)&:=\int_0^t \frac{X_0(s)[1-X_0(s)]}{N}ds,\quad \text{and}\\
\sigma(t)&:=\inf\{s>0, A(s)>t\}.
\end{align*}
Since
$$\int_0^{\sigma(t)}\frac{X_0(s)(1-X_0(s))}{N}ds=t,$$ we deduce that
\begin{align*}
\frac{d\sigma(t)}{dt}&=\frac{N}{{\tilde{X}_0}(t)(1-{\tilde{X}_0}(t))},\quad \text{provided we let}\\
{\tilde{X}_0(t)}&:=X_0({\sigma(t)}).
\end{align*}
Finally
$$\sigma(t)=\int_0^t\frac{N}{{\tilde{X}_0}(s)(1-{\tilde{X}_0}(s))}ds,$$
and if we let
$$\tilde{M}_1(t):=M_1(\sigma(t)),$$
we deduce from the above SDE for the process $X_0$ that 
$${\tilde{X}_0}(t)=x+N\int_0^t\frac{\alpha \tilde{M}_1(s)-\lambda}{1-{\tilde{X}_0}(s)}ds+B(t),$$
where $B(t)$ is a new standard Brownian motion (we use the same notation as above, which is a
slight abuse).

At the $k$--th step of our iterative procedure, $k\ge1$, we let ${\tilde{X}_0}$ start from $x-\sum_{j=1}^{k-1}\delta_j$, and we stop the process ${\tilde{X}_0}$ at the first time that it reaches the level  $x-\sum_{j=1}^{k}\delta_j$.
We will choose not only the sequence $\delta_k$, but also the sequence $s_k$ in such a way that
we can deduce from Lemma \ref{Bdescend}  (see \eqref{petit} and \eqref{trespetit} below)  that for each $1\le k\le K$ ($K$ to be defined below),
\begin{equation}\label{controlY}
\P\left(\inf_{0\le s\le s_k}\{\Theta_ks+B_s\}\le-\delta_k,  \sup_{0\le s\le s_k}\{\Theta_ks+B_s\}\le \tilde{\mu}\right)>
\frac{1}{3}.
\end{equation}
We shall make sure that
\begin{equation}\label{mintt}
\Theta_0=\beta', \text{ and }\Theta_k-\Theta_{k-1}\ge \tilde{\ep}+\lambda s'_k,
\end{equation}
with  $s'_k:=\sigma(s_k)$ and $\tilde{\ep}$ defined by \eqref{ep}, 
so that we deduce from Lemma \ref{M_1 pas trop rapide} and our choice of 
$\tilde{\ep}$ that
\begin{equation}\label{controlM}
\P(\sup_{0\le s\le s'_k}M_1(s)\le \Theta_k\Big| M_1(0)\le \Theta_{k-1})\ge3/4.
\end{equation}
The fact that with positive probability $X_0$ goes down to $\delta'$, while staying below $x+\tilde\mu$ and $M_1$ remaining under control,  will follow from a combination of \eqref{controlY} and \eqref{controlM}, provided
we show that we can choose the two sequences $\delta_k$ and $s_k$ for $k\ge1$ in such a way that
not only \eqref{controlY} holds, but also that there exists $K<\infty$ such that
$$x-\sum_{k=1}^K\delta_k\le\delta'.$$

 Since during the $k$--th step we are considering the event that $X_0(t)\le x+\tilde{\mu}$ i.e. $1-X_0(t)\ge a$, and also $X_0(t)\ge
 x-\sum_{j=1}^{k}\delta_j$,
 we have that 
 $$s'_k\le \frac{N}{a(x-\sum_{j=1}^{k}\delta_j)}s_k,$$
 so that we may, in accordance with \eqref{mintt}, make the following choice of $\Theta_k$ in terms of $\{\delta_j, s_j,\ 1\le j\le k\}$~:
 $$\Theta_k:=  \beta'+ k\tilde{\ep}+ N\frac{\lambda}{a}\sum_{j=1}^k  \frac{s_j}{x-\sum_{i=1}^{j}\delta_i}.$$
 We first want to ensure that (the reason for 0.4 will be made clear below)
 $$\frac{\delta_k}{\sqrt{s_k}}+\Theta_k\sqrt{s_k}\le 0.4,$$ which we achieve by requesting both that
 \begin{equation}\label{delta}
 \delta_k=0.2\sqrt{s_k}
 \end{equation} and
 \begin{equation}\label{t}
 \Theta_k\sqrt{s_k}\le 0.2\Leftrightarrow s_k\le\left(\frac{0.2}{\Theta_k}\right)^2.
 \end{equation}
On the other hand, we shall also request that for each $j\ge1$,
$$\frac{\delta_j}{x-\sum_1^j\delta_i}\le1\Leftrightarrow \delta_j\le\frac{1}{2}(x-\sum_{i=1}^{j-1}\delta_i).$$
This combined with \eqref{delta} implies that
$$\frac{s_j}{x-\sum_{i=1}^j\delta_i}\le 25 \delta_j.$$
Consequently
\begin{align*}
\beta' \le \Theta_k&\le \beta'+ k\tilde{\ep}+25N\frac{\lambda }{a}\left(\sup_{1\le j\le k}\delta_j\right) k .
\end{align*}
 Moreover, a combination of \eqref{delta} and \eqref{t} yields
\begin{align*}
\delta_j &= 0.2 \sqrt{s_j} \leq \frac{(0.2)^2}{\Theta_j} \\
& \leq (25 \beta')^{-1},
\end{align*}
and from the above inequality follows
\begin{align}\label{tetak}
\Theta_k&\le \beta'+D_N k,\\
\text{with }D_N&= \tilde{\ep}+\frac{N\lambda}{a\beta'}.\nonumber
\end{align}

Finally this leads us to choose
\begin{align}\label{kappa}
\delta_k&=\inf\left(\frac{\kappa}{(\beta'+D_N k)},\frac{1}{2}(x-\sum_{j=1}^{k-1}\delta_j)\right),\\
s_k&= 25 \delta_k^2.\nonumber
\end{align}
It still remains to choose $\kappa$, which will be done below. Note that
\eqref{delta} + \eqref{t} request us to make sure that
$\kappa \le \frac{1}{25} $.

We now have 
\begin{lemma}
$\exists K>0,$ $ \forall k>K,$
$$\delta_k=\frac{1}{2}(x-\sum_{j=1}^{k-1}\delta_j).$$
\end{lemma}

\bpf
We first show that for $k\ge2$,
$$\frac{1}{2}(x-\sum_{j=1}^{k-1}\delta_j)\le\frac{\kappa}{(\beta'+D_N k)}\Rightarrow
\frac{1}{2}(x-\sum_{j=1}^{k}\delta_j)<\frac{\kappa}{(\beta'+D_N (k+1))}.$$
Indeed, if the above left inequality holds, then
$$\frac{\frac{\kappa}{(\beta'+D_N (k+1))}} {\frac{1}{2}(x-\sum_{j=1}^{k-1}\delta_j)} \geq \frac{\frac{\kappa}{(\beta'+D_N (k+1))}} {\frac{\kappa}{(\beta'+D_N k)}}  
> \frac{1}{2},$$
where the last inequality follows easily from $k\ge2$. Consequently
\begin{align*}
\frac{\kappa}{(\beta'+D_N (k+1))}&>\frac{1}{2}\delta_k\\
&=\frac{1}{4}\left(x-\sum_{j=1}^{k-1}\delta_j\right)\\
&=\frac{1}{2}\left(x-\sum_{j=1}^{k}\delta_j\right).
\end{align*}
Finally there exists $K'\ge1$ such that $$x-\sum_{j=1}^{K'}\frac{\kappa}{(\beta'+D_N j)}<0.$$
Therefore for some $k\le K'$,
$$\frac{\kappa}{(\beta'+D_N k)}>\frac{1}{2}\left(x-\sum_{j=1}^{k}\delta_j\right).$$
\epf

This means that at each $k>K$, ${\tilde{X}_0}$ progresses by a step equal to half the remaining distance to zero. 
Consequently $\exists c>0$ such that $x_k=x-\sum_{j=1}^{k}\delta_j \leq c2^{-k}$. We are looking for the smallest integer $\overline{k}$ such that
$c2^{-\overline{k}}\le\delta'$, $\delta'$ to be specified below, which implies that 
$$\overline{k}-1 < \frac{\log(c)-\log(\delta')}{\log(2)}\le\overline{k}.$$
Consequently, since we may as well assume that $\delta'\le1/2$,
\begin{align*}
\overline{k}&\le 1+\frac{\log(c)}{\log(2)}+\frac{\log(1/\delta')}{\log(2)}\\
&\le \left[\log(2)\right]^{-1}\left(2+\frac{\log(c)}{\log(2)}\right)\log\left(\frac{1}{\delta'}\right).
\end{align*}
Combining this estimate with \eqref{tetak}, we deduce that there exists a constant $D'_N$ such that
$$\Theta_{\overline{k}}\le \beta'+D'_N \log\left(\frac{1}{\delta'}\right).$$
 Hence 
 there exists a $\delta' \leq \delta\wedge1/2 $ such that $\delta'\Theta_{\overline{k}}\le\eps$.
 If we now check that the probability of the previous path is bounded below by a positive constant,
 we will have that with a positive constant, at the end of the $\overline{k}$--th step, both $X_0\le\delta'\le\delta$, and $M_1\le\Theta_{\overline{k}}$, hence $X_0M_1\le\eps$, which puts us in a position to apply Proposition \ref{clef}.

Given the choice that we have made for $\tilde{\ep}$, see \eqref{ep}, it suffices to make sure that
\begin{equation}\label{petit}
\sqrt{\frac{2}{\pi}}\left(\frac{\delta_k}{\sqrt{s_k}}+\Theta_k\sqrt{s_k}\right)<1/3,\ \forall k\ge1,
\end{equation}
as well as
\begin{equation}\label{trespetit}
2\exp\left[-\frac{1}{2}\left(\frac{\tilde{\mu}}{\sqrt{s_k}}-\Theta_k\sqrt{s_k}\right)^2\right]<1/3,\ \forall k\ge1.
\end{equation}
Since $3^{-1}\sqrt{\pi/2}>0.4$, \eqref{delta}+\eqref{t} implies \eqref{petit}.

On the other hand, (\ref{trespetit}) is equivalent to
\begin{equation}\label{ineq}
\left( \frac{\tilde{\mu}}{\sqrt{s_k}}-\Theta_k\sqrt{s_k} \right)^2 > 2 \log 6.
\end{equation}
But we have
\begin{lemma}
A sufficient condition for \eqref{ineq} is that
\begin{equation}\label{ineq2}
\kappa\le1\wedge\frac{\beta'\tilde{\mu}}{25+10\log6}.
\end{equation}
\end{lemma}
 \bpf
 It follows from \eqref{ineq2}
 \begin{align*}
 \tilde{\mu}(\beta'+D_Nk)&>(25+10\log6)\kappa\\
 &\ge\left(\frac{25\Theta_k}{\beta'+D_Nk}+10\log6\right)\kappa\\
 &\ge\frac{25\Theta_k}{\beta'+D_Nk}\kappa^2+10(\log6)\kappa,\\
 \frac{\tilde{\mu}}{5\kappa/(\beta'+D_Nk)}&>\Theta_k\frac{5\kappa}{\beta'+D_Nk}+2\log6.
 \end{align*}
 Finally \eqref{ineq} follows from the last inequality, \eqref{kappa} and \eqref{delta}.  
 \epf
 
 We therefore choose
 $$\kappa = \frac{1}{25} \wedge\frac{\beta'\tilde{\mu}}{25+10\log6}.$$
 
 We can now conclude that
 \begin{proposition}\label{step2}
 Suppose that $X_0(0)\le x_{\max}$ and $M_1(0)\le\beta'$. Let
 $$T_{\delta'}=\inf\{s>0,\ X_0(s)\le\delta'\}.$$
 Then
 $$\P\left(T_{\delta'}\le t_2,\ X_0(T_{\delta'})\times M_1(T_{\delta'})\le\ep\right)\ge\left( \frac{1}{12} \right)^{\overline{k}_{\max}}:=p_{trans} ,$$
 with $t_2=25\overline{k}_{\max}$, and $\overline{k}_{\max}$ is the number of steps needed to reach $\delta'$ in the above procedure, while starting from $x_{\max}$.
 \end{proposition}
\bpf
It follows from \eqref{controlM}, \eqref{petit}, \eqref{trespetit}, Lemma \ref{Bdescend} and again Lemma \ref{trivial} that the $k$--th step in the above procedure happens with probability at least $1/12$. It remains to exploit the Markov property, like at the end of the proof of Lemma \ref{Mini}. \epf

  \subsection{Conclusion}
  
  Proposition \ref{stepun} shows that, if we start with $M_1<\beta$, with probability
  $p_{init}$ we need to wait at most a length of time $t_1$ for the pair $(X_0,M_1)$
  to reach the set $[0,x_{max}]\times[0,\beta']$. Proposition \ref{step2} shows that starting from that set, with probability $p_{trans}$ we need to wait at most a length of time $t_2$
  for $(X_0,X_0M_1)$ to reach the set $[0, \delta']\times[0,\epsilon]$, with $\delta'<\delta$.
  But from Proposition \ref{clef}, starting from this last set, we have a probability $p_{fin}$ to reach $0$ 
  during an interval of time of length ${t_3}$.
  
  So to sum up, using again the strong Markov property of the system, we have
  
  \begin{proposition}\label{iteration}
  For any finite stopping time $T$, if $M_1(T) \leq \beta$,  then
  $$\P( T_0 < T+ t_1+t_2+{t_3}) \geq p_{fin} p_{trans} p_{init} >0. $$
  \end{proposition}
 
  Moreover Lemma \ref{Mini} implies that this situation will happen infinitely many times as long as the ratchet does not click, which implies Theorem \ref{th}, exploiting again the strong Markov property of the solution of \eqref{un}.

\section{Proof of Theorem \ref{thdeu}}

This final section is devoted to the proof of 
Theorem \ref{thdeu}.

 We first note that the reasoning of section 5 can be done with any initial value $\rho$ for $M_1$, instead of $\beta$. That is to say, with $S_\rho^t=\inf\left\{s>t,M_1(s) \le \rho \right\} $ (and $S_\rho=S_\rho^0$),

  \begin{lemma}\label{le5.1}
   $\exists$ $t^\rho_1,t^\rho_2,t^\rho_3<\infty$, and $p^\rho_{init},p^\rho_{trans},p^\rho_{fin}$ $>0$ such that
  
  $$\P( T_0 < S_\rho^t+ t^\rho_1+t^\rho_2+{t^\rho_3}) \geq p^\rho_{init} p^\rho_{trans} p^\rho_{fin}\text{ .}$$
  
  \end{lemma}
Choosing $\rho= \frac{\ep}{\delta} \vee \frac{2 \lambda}{\alpha}$, we have~:

\begin{lemma}\label{Smfini}
There exist $K$, $\tilde{p}>0,$ such that for any initial condition in the set $\x_{\delta}$,
$$\P(T_0\wedge S_\rho  \le K)\ge\tilde{p}.$$

\end{lemma}

\bpf
We are going to argue like in the proof of Lemma \ref{Mini}.
We introduce the process $\lbrace Y_{s}, s\ge 0 \rbrace$, which is the solution of the following system :

\begin{equation}
\left\{
\begin{aligned}
 dY_{s} &=  \frac{\alpha \ep}{2} ds + \sqrt{\frac{Y_{s}(1-Y_{s})}{N}}dB_0(s), \\
 Y_0 &= 0.
\end{aligned}
\right.
\end{equation}

For any $0 \le u \le 1,$ let $$R_u= \inf\{s \geq 0, Y_s=u \}.$$
Since $\frac{\alpha \ep}{2} >0$ we deduce that there exist $L>0$, $p>0$ such that $ \P(R_1 \leq L ) \geq p >0$. We choose $K=L + {t_3}$, where ${t_3}$ has been defined in Proposition \ref{clef}.
 
 Now there are several possibilities :
\begin{description}
\item{Case 1.} $\inf_{ 0 \leq s \le L} M_1(s) \leq \rho$, then $S_\rho < L<K$.
\item{Case 2a.} $\inf_{ 0 \leq s \le L} M_1(s) \geq \rho$ and $\inf_{ 0 \leq s \le L} X_0(s) M_1(s) \leq \ep $.  Then there exists $t < L$ such that $X_0(t) M_1(t) \leq \ep$ (which implies $X_0(t) \le \delta$, because $ M_1(t) \geq \rho \ge \frac{\ep}{\delta}$). In that case we can use Proposition \ref{clef}, and we have $\P(T_0 \le K)\ge p_{fin} >0$, which implies $\P(T_0 \wedge S_\rho \le K)
\ge p_{fin} >0$.
\item{Case 2b.} $\inf_{ 0 \leq s \le L} M_1(s) \geq \rho$ and $\inf_{ 0 \leq s \le L} X_0(s) M_1(s) \geq \ep $. In that last case we have (using first $X_0 \ge \frac{\ep}{M_1}$ combined with  $\alpha M_1- \lambda \ge \lambda >0$, and next $-\frac{\lambda}{M_1(s)}\ge-\frac{\alpha}{2}$)
\begin{align*}
\inf_{ 0 \leq s \le L}(\alpha M_1(s) - \lambda )X_0(s) &\geq \inf_{ 0 \leq s \le L} \ep(\alpha - \frac{\lambda}{M_1(s)}) \\
& \geq \frac{\alpha \ep}{2},
\end{align*}
and consequently we can use the comparison theorem (Lemma \ref{comparaison}), which implies that $\forall s \in \left[0,L\right]$, $X_0(s) \ge Y_s$. Then $\P(T_1 \le L) \ge p >0$. But when $X_0$ hits $1$, $M_1$ hits $0$. Hence $\P(S_{\rho} \le L) \ge p >0$. 
\end{description}

We may now conclude that there exists $\tilde{p}>0$ such that $$\P(T_0 \wedge S_\rho \le K ) \ge \tilde{p}.$$

\epf

We deduce from the two above Lemmas :

\begin{corollary}\label{cor5.2}
  There exists $\overline{K}<\infty$, and $\overline{p}>0$ such that, for any initial condition in $\x_\delta$ for some $\delta>0$,
  $$\P( T_0 \le \overline{K}) \geq \overline{p}.$$
\end{corollary}

We can now proceed with the

\noindent{\sc Proof of Theorem \ref{thdeu}}
We deduce from Corollary \ref{cor5.2} and the strong Markov property that for all $n\ge0$, 
$\P(T_0>n\overline{K})\le(1-\overline{p})^n$. Consequently
\begin{align*}
\E[e^{\rho T_0}]&\le \sum_{n=0}^{\infty} e^{(n+1)\rho\overline{K}} 
\P(n\overline{K}\le T_0\le(n+1)\overline{K})\\
&\le   \sum_{n=0}^{\infty} e^{(n+1)\rho\overline{K}}(1-\overline{p})^n  \\
&= e^{\rho\overline{K}} \sum_{n=0}^{\infty} \left(e^{\rho\overline{K}}(1-\overline{p})\right)^n\\
&<\infty,
\end{align*}
provided $\log(1-\overline{p})+\rho\overline{K}<0$, in other words 
$\rho<\overline{\rho}:=-\log(1-\overline{p})/\overline{K}$.
\epf

\paragraph{Acknowledgements} The authors thank Jean--St\'ephane Dhersin, Peter Pfaffelhuber
and Anton Wakolbinger for valuable discussions concerning this work, as well as anonymous Referees, whose reports helped us to improve the exposition.

\end{document}